\documentclass[11pt,a4paper,enabledeprecatedfontcommands]{scrartcl}

\usepackage[utf8]{inputenc}
\usepackage[T1]{fontenc}
\usepackage[english]{babel}

\usepackage{xcolor}
\usepackage[bookmarks,
bookmarksnumbered,
linktocpage,
colorlinks=true,
citecolor=blue,
linkcolor=violet
]{hyperref}

\usepackage[pdftex]{graphicx}
\usepackage{latexsym}
\usepackage{amsfonts}
\usepackage{setspace}
\usepackage{amsmath,amsthm,amssymb,comment}
\allowdisplaybreaks
\usepackage{geometry}

\newtheorem{theo}{Theorem}[section]
\newtheorem{lem}[theo]{Lemma}

\newtheorem{defi}{Definition}[section]	   
                  
\numberwithin{equation}{section}

\newcommand{\R}{\mathbb{R}} 
\newcommand{\eps}{\varepsilon} 
\newcommand{\f}[2]{\frac{#1}{#2}} 

\newcommand{\be}{\begin{equation} \label}
\newcommand{\ee}{\end{equation}}
\newcommand{\bea}{\begin{eqnarray}\label}
\newcommand{\eea}{\end{eqnarray}}
\newcommand{\bas}{\begin{eqnarray*}}
\newcommand{\eas}{\end{eqnarray*}}
\newcommand{\bit}{\begin{itemize}}
\newcommand{\eit}{\end{itemize}}
\newcommand{\nn}{\nonumber}
\newcommand{\lbal}{\left\{ \begin{array}{l}}
\newcommand{\lball}{\left\{ \begin{array}{ll}}
\newcommand{\ear}{\end{array} \right.}

\newcommand{\pO}{\partial\Omega} 

\newcommand{\Proof}{\vspace*{2mm} {\sc Proof.} \quad} 
\newcommand{\abs}{\\[5pt]} 

\newcommand{\uy}{\underline{y}} 
\newcommand{\Om}{\Omega} 
\newcommand{\io}{\int\limits_\Omega} 
\newcommand{\bom}{\overline{\Omega}} 
\newcommand{\uw}{\underline{w}} 
\newcommand{\ow}{\overline{w}} 
\newcommand{\norm}[1]{\left\lVert#1\right\rVert} 

\begin{document}

\enlargethispage{10mm}
\title{Solutions to a chemotaxis system with spatially heterogeneous diffusion sensitivity}
\author{
Gregor Flüchter\footnote{gmf@mail.upb.de}\\
{\small Institut f\"ur Mathematik, Universit\"at Paderborn,}\\
{\small 33098 Paderborn, Germany} }
\date{}
\maketitle
\begin{abstract}
\noindent 
  We consider a parabolic-elliptic Keller-Segel system with spatially dependent diffusion sensitivity
  \begin{eqnarray*}
	\left\{ \begin{array}{l}
	u_t = \nabla \cdot (|x|^\beta \nabla u) - \nabla \cdot (u\nabla v), \\[1mm]
	0 = \Delta v - \mu + u,
	\qquad \mu:=\frac{1}{|\Omega|} \int\limits_\Omega u,
 	\end{array} \right.
	\qquad \qquad (\star)
  \end{eqnarray*}
  under homogeneous Neumann boundary conditions in the ball $\Omega=B_R(0)\subset \mathbb R^n$.\\
For $\beta>0$ and radially symmetric Hölder continuous initial data, we prove that there exists a pointwise classical solution to $(\star)$ in $(\Omega\setminus \{0\})\times (0,T)$ for some $T>0$.\\
For radially decreasing initial data satisfying certain compatibility criteria, this solution is bounded and unique in $(\Omega\setminus \{0\})\times (0,T^*)$ for some $T^*>0$.\\
Moreover, for $n \geq 2$ and sufficiently accumulated initial data, there exists no solution $(u,v)$ to $(\star)$ in the sense specified above which is globally bounded in time.\\[5pt]

\end{abstract}
\newpage

  \pagestyle{headings}

\section{Introduction}

In microbiological processes, it is common for organisms to interact with their environment via positive chemotaxis, that is the tendency to move in the direction of the gradient of some signal substance. This behavior has been documented as early as 1881 for Bacterium termo and Spirillum tenue moving toward oxygen-producing plant cells \cite{Engelmann}.\\
In order to make such biological systems accessible to quantitative analysis and outline the governing factors of structural evolution, at the beginning of the 1970s Keller and Segel proposed a system of the form
\be{KSpara}
	\left\{ \begin{array}{l}
	u_t = \nabla \cdot (D(u,v) \nabla u) - \nabla \cdot (S(u,v)\nabla v), \\
	v_t =  \Delta v + u - v,
	\end{array} \right. 
\ee
see \cite{KelSe1} and \cite{KelSe2}, to describe the behavior of slime mold aggregation and Escherichia coli as outlined in \cite{Adler}, respectively. Herein, the bacteria concentration $u=u(x,t)$ is subject not only to chemotaxis toward a self-produced attractant with density $v=v(x,t)$ but also to diffusion in the form of Brownian motion. The expressions $D(u,v)$ and $S(u,v)$ represent the diffusive and chemotactic sensitivity, respectively.\\
This system and variations thereof have various applications \cite{HilPa}, including pattern formation in bacterial colonies \cite{wtmmbb}, tumor invasion processes \cite{ChaLo} and embryogenesis \cite{pmo}.\\
In line with experimental observations, even the prototypical setting with $D(u,v)\equiv 1$ and $S(u,v)\equiv u$ considered in bounded domains with no-flux boundary conditions has been shown to exhibit aggregation phenomena already, in their most extreme form represented by finite-time blow-up. While solutions always exist globally in time and are bounded in the spatially one-dimensional case \cite{OsaYa}, for two-dimensional balls, a critical mass phenomenon arises: For all sufficiently regular and radially symmetric initial data $u_0$ with total mass $\io u_0 < 8\pi$, solutions are global and bounded \cite{nsy}, whereas if $\io u_0 > 8\pi$, finite-time blow-up is possible (\cite{HerVe}, \cite{MizWi}). For balls as domains in higher dimensions, initial data with arbitrary total mass leading to blow-up can be constructed \cite{Win2013}.\\
On account of the fact that in numerous biological applications, the chemo-attractant dissipates much faster than the microbes move, by \cite{JagLu} we may consider the parabolic-elliptic system given by
\be{KSparaell}
	\left\{ \begin{array}{ll}
	u_t = \nabla \cdot (D(u,v) \nabla u) - \nabla \cdot (S(u,v)\nabla v), \quad &x\in\Om,~t>0,\\
	0 =  \Delta v + u - \mu, \qquad \mu=\f1{|\Om|}\io u, \quad &x\in\Om,~t>0,
	\end{array} \right. 
\ee
in a bounded domain $\Om\subset \R^n$ as a relevant limit case of (\ref{KSpara}).
Concerning global boundedness and finite-time blow-up, the results are similar to the fully parabolic system (\ref{KSpara}) \cite{nagaisen}.
The matter of deducing local existence in both systems in absence of degeneracies, particularly for constant and linear sensitivities $D$ and $S$, respectively, has for example been considered in \cite{biler98}. Therein, heat semigroup theory is employed in order to obtain a mild solution which can be shown to have nice regularity properties in the interior (\cite[III.12]{LSU}) in turn. By means of a fixed point argument, such results have been extended to possibly degenerate cases for sufficiently regular $D>0$ and $S\geq 0$ dependent on $u$ (\cite{dw2010}). Related systems for nutrient taxis with $D(u,v)\equiv uv$ and $S(u,v)=\Psi(u)v$ with $\Psi$ asymptotically growing quadratically at most as well as such with $D(u)$ and $S(x,u,v)$ supposed to be not too singular in certain ways have been investigated \cite{WinLinf} and \cite{Winlogsing}, respectively, and, utilizing approximations, have been shown to possess global weak solutions. Further variants of (\ref{KSpara}) and (\ref{KSparaell}) containing density-- and signal--dependent diffusion degeneracies have been discussed (see, e.g., \cite{tato23}, \cite{ChiMiTo20}, \cite{MizOnYo}, \cite{Jiang18}, \cite{fujiang21}, \cite{Iyo22}, \cite{Iyo23}).\\
In this manuscript, we shall consider a spatially dependent diffusion sensitivity generalized by the prototype $D(x)\equiv |x|^\beta$, $x \in \Om$, for $\beta>0$.
This leads to the problem
\be{DS}
    	\left\{ \begin{array}{rcll}
	u_t &=& \nabla \cdot (|x|^\beta \nabla u) - \nabla \cdot (u\nabla v), 
	\qquad & x\in\Omega, \ t>0, \\[1mm]
	0 &=& \Delta v - \mu + u, \quad \mu:=\frac{1}{|\Omega|} \io u,
	\qquad & x\in\Omega, \ t>0, \\[1mm]
	& & \hspace*{-10mm}
	\frac{\partial u}{\partial\nu}=\frac{\partial v}{\partial\nu}=0,
	\qquad & x\in\pO, \ t>0, \\[1mm]
	& & \hspace*{-10mm}
	u(x,0)=u_0(x),
	\qquad & x\in\Omega, 
 	\end{array} \right.
\ee
posed in $\Omega=B_R(0)\subset\R^n$, $n\ge 1$, $R>0$,
fulfilling
 \be{init}
	u_0\in C^0_{rad}(\bom) := \Big\{ \varphi\in C^0(\bom) \ \Big| \ \varphi \mbox{ is radially symmetric } 
	\Big\}
	\quad \mbox{is nonnegative with $u_0\not\equiv 0$}.
\ee
We shall remark that in the contexts we consider, actually $\mu=\frac{1}{|\Omega|} \io u_0$ will be proven to hold, so that $\mu$ may be viewed as a constant.\\
This system can be interpreted as a prototype for describing biological applications where the motility of a cell or bacteria population is impaired near the origin. For instance, we find this to be the case when coagulation mechanisms are present. Their significance not just for structural healing but also in the context of immune responses  for invertebrates has among others been established in \cite{osabiol} and \cite{Wangbiol}. In mammals, the coagulation system was long thought to be important exclusively for haemostasis. However, nowadays it is commonly recognized that coagulation contributes to the effective elimination of bacteria in those organisms as well \cite{berends}. In fact, besides restricting the motility of bacteria, coagulation triggers the release of bradykinin which interacts with macrophages to emit chemo-attractants supporting the immune response \cite{clotthum}. Furthermore, fibrinogen releases fibrinopeptides, chemo-attractants to aid clotting \cite{skogen}. Thus in this example already, there are multiple chemotactic dynamics at play wherein heterogeneous environments roughly as described in (\ref{DS}) might occur.\\
The mathematical analysis of (\ref{DS}) however is accompanied by notable difficulties. Calculating
\[
\nabla \cdot (|x|^\beta \nabla u)=|x|^\beta \Delta u + (\nabla |x|^\beta)\cdot \nabla u
\]
reveals that for one, we are dealing with a spatial diffusion degeneracy which in Keller-Segel type systems appears to be without precedent in literature, and moreover, at least for $\beta<1$ singular behavior of $(\nabla |x|^\beta)\cdot \nabla u$ at $x=0$ is to be expected. This already indicates that at least generally, we should not assume to be able to obtain a classical solution of (\ref{DS}) in $\bom \times (0,T)$ for some $T>0$; instead, we either have to resort to weak solution concepts or at least omit the spatial point $x=0$. Our results feature the latter.\\

First we formulate a statement on local existence of classical solutions in $(\bom\setminus \{0\})\times(0,T_0)$ satisfying mass conservation.

\begin{theo}\label{maintheo1}
Let $n\geq 1$, $R>0$, $\Om=B_R(0)\subset \R^n$ and $\beta>0$, and write $\Om_0:=\bom\setminus \{0\}$.\\
Then for $\theta\in (0,1)$ and $u_0 \in C^\theta(\bom)$ complying with (\ref{init}), there exists a radially symmetric classical solution (u,v) of (\ref{DS}) in the sense of Definition \ref{classsol}, fulfilling
\be{uvreg}
\begin{cases}
u \in C^0(\Om_0 \times [0,T_0))\cap C^{2,1}(\Om_0 \times (0,T_0))\\
v \in C^{2,0}(\Om_0 \times (0,T_0))
\end{cases}
\ee
for some $T_0\in(0,\infty]$. This solution has the properties that $u$ is nonnegative and satisfies the mass conservation property, that is
\be{massctheo}
\io u(\cdot,t)=\io u_0=:m \qquad \text{for all}\quad t\in(0,T_0).
\ee
\end{theo}
\vspace{10mm}
The second theorem includes a result on local boundedness and uniqueness. In order to accomplish this, we need to impose much stronger requirements on the initial data.

\begin{theo}\label{maintheo2}
Suppose that the conditions of Theorem \ref{maintheo1} hold, and let $(u,v)$ denote the classical solution to (\ref{DS}) established therein.\\
Assume that additionally $u_0 \in C^{1+\theta}(\bom)$ for some $\theta \in (0,1)$ is radially decreasing and has the properties that
\be{Omboundary}
u_0=0 \qquad \text{and} \qquad \nabla u_0 \cdot \nu =0 \qquad \text{on} \quad \pO
\ee
as well as
\be{xlimit}
|\nabla u_0(x)| \leq C_0 |x|^{n-1+\theta}
\ee
for some $C_0>0$.\\
Then for $T_0>0$ as in Theorem \ref{maintheo1} we have that 
\be{radiallydec}
u\in C^{1,0}(\Om_0\times[0,T_0)) \quad \text{is radially decreasing.}
\ee
Moreover, for some $T^*\in(0,T_0]$ and each $T\in(0,T^*)$, there exists $C=C(T)>0$ such that
\be{realuboundtheo}
u(x,t)\leq C \qquad \text{for all} \quad (x,t)\in\Om_0\times[0,T].
\ee
If additionally $n \geq 2$, there exists a unique solution $(u,v)$ of (\ref{DS}) in $\Om_0\times[0,T^*)$ fulfilling
\be{regularityclasses}
\begin{cases}
u \in C^0(\Om_0 \times [0,T^*))\cap C^{2,1}(\Om_0 \times (0,T^*)),\\
v \in C^{2,0}(\Om_0 \times (0,T^*)),
\end{cases}
\ee
which has the properties that
\be{unicondtheo}
\io v(\cdot,t) = 0, \quad 0\leq u \in L^\infty(\Om\times (0,T)) \quad \text{and} \quad \nabla v \in L^\infty(\Om\times(0,T);\R^n)
\ee
for all $T\in(0,T^*)$.
\end{theo}
\vspace{5mm}

We close with a result ruling out global bounded solutions for initial mass distributions concentrated adequately close to the origin.

\begin{theo}\label{maintheo3}
Let $n\geq 2$, $R>0$, $\Om=B_R(0)\subset \R^n$ as well as $\beta>0$, and assume that $u_0$ complies with (\ref{init}).\\
Then for $m:=\io u_0$ and each $m_0\in(0,m]$, there exists $r_0=r_0(m_0,m,R,\beta)>0$ such that if
\be{initmasscon}
\int_{B_{r_0}(0)}u_0 \geq m_0,
\ee
there is no global classical solution $(u,v)$ of (\ref{DS}) fulfilling
\be{uvreginfty}
\begin{cases}
u \in C^0(\Om_0 \times [0,\infty))\cap C^{2,1}(\Om_0 \times (0,\infty))\\
v \in C^{2,0}(\Om_0 \times (0,\infty))
\end{cases}
\ee
such that for each $T\in(0,\infty)$
\be{allisbound}
0\leq u \in L^\infty(\Om\times(0,T)) \quad \text{and} \quad \nabla v \in L^\infty(\Om\times(0,T);\R^n). 
\ee
\end{theo}
\vspace{5mm}

{\bf Outline of arguments.} 
The main idea is to transform the Keller-Segel type Neumann boundary value problem to a Dirichlet problem for which we are able to obtain a local solution, see subsections \ref{traforward} and \ref{auxproblem}, and then retransform in order to acquire a solution of (\ref{DS}) (Subsection \ref{re}).\\
In Section \ref{bu}, we are concerned with ruling out global boundedness for sufficiently large $\beta>0$ and properly concentrated initial data. The main idea here is to attach singular weights to the mass accumulation function $w$ and thus construct a generalized moment functional. This functional is bounded, but shown to explode in finite time under the assumption that $w$ solves (\ref{0w}) globally with $w_s$ bounded locally time, implying that this cannot be the case.\\
The main theorems can then be obtained mainly by collecting previous results. One needs to be cautious how the conditions imposed on the initial data of the original system (\ref{DS}) translate to those in (\ref{0w}) though.

 \newpage  

\section{Existence of solutions}

We shall establish the existence of sufficiently smooth solutions to (\ref{DS}).\\

In the scenario at hand, that represents a particular challenge. Not only is the term $|x|^\beta$ not differentiable at $0$ for $0<\beta\leq 1$, but moreover the coefficient of the Laplacian of $u$ vanishes at $x=0$ for all $\beta>0$, implying a diffusion degeneracy. Whereas examples of possible degeneracies depending on $u$ or even on $u$ and $v$ have for instance been discussed in \cite{dw2010}, \cite{WinSin} and \cite{WinLinf}, to the author's knowledge no case of a spatially dependent diffusion degeneracy in such systems with Neumann boundary conditions has yet been addressed in standard literature.\\

Our approach in principle relies on a strategy usually employed to detect blow-up in parabolic-elliptic Keller-Segel type chemotaxis systems, introduced by Jäger and Luckhaus in \cite{JagLu}. We remark that approaches based on an analysis of mass accumulation functions can be found in numerous works on related problems (confer \cite{MaoLi}, \cite{tumuzh}, \cite{Wangli}, \cite{duliu23}, \cite{liuli21}). In most precedent cases of this type, however, considerations in this regard concentrate on the construction of exploding solutions, with only few exceptions (confer \cite{biler08}); we emphasize that in the present manuscript already the mere construction of solutions operates on a level of cumulated densities. Therefore, it is crucial to us that not only radial symmetry but also mass is conserved for adequately regular solutions.
\abs

\subsection{Transforming with respect to a weaker solution concept}\label{traforward}

Let us first define our solution concept.

\begin{defi}\label{classsol}
Suppose that $n\geq 1$, $R>0$, $\Om=B_R(0) \subset \R^n$, $T>0$, $\Om_0:=\bom \setminus \{0\}$, and let $u_0$ comply with (\ref{init}). We call a pair of functions $(u,v)$ satisfying
\be{C21A}
\left\{ \begin{array}{l}
	u \in C^0(\Om_0 \times [0,T))\cap C^{2,1}(\Om_0 \times (0,T)),\\[1mm]
	v \in C^{2,0}(\Om_0 \times (0,T)),
	\end{array} \right.
\ee
and solving (\ref{DS}) pointwise in $(\bom \setminus \{0\})\times [0,T)$ a classical solution of (\ref{DS}) in $(\bom \setminus \{0\}) \times [0,T)$.\\
Analogously, we call $(u,v)$ with the property (\ref{C21R}) solving (\ref{DSR}) pointwise in $(0,R]\times [0,T)$ a classical solution of (\ref{DSR}) in $(0,R] \times [0,T)$.\\
\end{defi}

For radially symmetric solution of the system (\ref{DS}), we rewrite (\ref{DS}) in radial coordinates. By writing $r:=|x|$, the pair of functions $(u,v)=(u(r,t),v(r,t))$ with
\be{C21R}
\left\{ \begin{array}{l}
	u \in C^0((0,R] \times [0,T))\cap C^{2,1}((0,R] \times (0,T)),\\[1mm]
	v \in C^{2,0}((0,R] \times (0,T)),
	\end{array} \right.
\ee
fulfills
\be{DSR}
    	\left\{ \begin{array}{rcll}
	u_t &=& \f1{r^{n-1}}(r^{n-1+\beta}u_r)_r - \f1{r^{n-1}}(r^{n-1} uv_r)_r, 
	\qquad & r\in(0,R), \ t>0, \\[1mm]
	0 &=& \f1{r^{n-1}}(r^{n-1}v_r)_r - \mu + u,
	\qquad & r\in(0,R), \ t>0, \\[1mm]
	& & \hspace*{-10mm}
	u_r=v_r=0,
	\qquad & r=R, \ t>0, \\[1mm]
	& & \hspace*{-10mm}
	u(r,0)=u_0(r),
	\qquad & r\in(0,R),
 	\end{array} \right.
\ee
pointwise in $(0,R]\times [0,T)$.\\
Analogously, we name $(u,v)$ with such properties a classical solution of (\ref{DSR}) in $(0,R] \times [0,T)$.\\

From here on further, we shall denote $r:=|x|$ and without risk of confusion write $(u,v)=(u(r,t),v(r,t))$ in the context of (\ref{DSR}).\\

In order to deal with these solutions defined on a non-compact space, we shall focus on bounded solutions.

\begin{lem}\label{vrlem}
Suppose that $n\geq 1,~R>0$, and let $u_0$ comply with (\ref{init}). Let $(u,v)$ with $u\in L^\infty((0,R)\times (0,T_0))$ be a classical solution of (\ref{DSR}) in $(0,R] \times [0,T_0)$.\\
If then
\begin{equation}\label{vr}
v_r(r,t)=\f1{r^{n-1}}\bigg(\f{\mu r^n}n-\int\limits_0^r \rho^{n-1}u(\rho,t)d\rho\bigg) \quad \text{for all }(r,t)\in (0,R] \times [0,T_0), 
\end{equation}
we have 
\begin{equation}\label{vrbound}
|v_r(r,t)|\leq Cr \qquad \text{for all }(r,t)\in (0,R] \times [0,T_0)
\end{equation}
with $C:=\f{2}n \cdot \norm{u}_{L^\infty((0,R)\times(0,T_0))}$, and hence in particular $v_r \in L^\infty((0,R)\times (0,T_0))$.\\
Moreover, for $n\geq 2$ the converse statement also holds true:\\
If $n\geq 2$ and $v_r \in L^\infty((0,R)\times (0,T_0))$, then necessarily (\ref{vr}).
\end{lem}

\Proof
Observe that $v_r$ defined as in (\ref{vr}) indeed complies with (\ref{DSR}) since
\[
(r^{n-1}v_r)_r=\bigg(\f{\mu r^n}n-\int\limits_0^r \rho^{n-1}u(\rho,t)d\rho\bigg)_r=\mu r^{n-1} - r^{n-1}u
\]
and thus
\[
0=\f1{r^{n-1}}(r^{n-1}v_r)_r - \mu + u,
\]
as well as
\[
v_r(R,t)=\f1{R^{n-1}}\bigg(\f{\mu R^n}n-\int\limits_0^R \rho^{n-1}u(\rho,t)d\rho\bigg)=0
\]
due to
\[
\mu=\frac{1}{|B_R(0)|} \int\limits_{B_R(0)} u = \f n{\omega_n R^n}\omega_n \int_0^R \rho^{n-1}u(\rho,t) d\rho = \f n{R^n}\int_0^R \rho^{n-1}u(\rho,t) d\rho.
\]
As a consequence of the boundedness of $u$ in $(0,R] \times [0,T_0)$, we now obtain that for $(r,t)\in (0,R]\times(0,T_0)$
\begin{align*}
|v_r(r,t)|&=\f1{r^{n-1}}\bigg|\f{\mu r^n}n-\int\limits_0^r \rho^{n-1}u(\rho,t)d\rho\bigg|\\
&=\f1{r^{n-1}}\bigg|\int\limits_0^r \rho^{n-1}(\mu-u(\rho,t))d\rho\bigg|\\
& \leq \f1{r^{n-1}}\int\limits_0^r \rho^{n-1}\norm{\mu-u}_{L^\infty((0,R)\times(0,T_0))}d\rho\\
& \leq Cr
\end{align*}
with $C:=\f{2}n \cdot \norm{u}_{L^\infty((0,R)\times(0,T_0))}$, and therefore 
\[
\norm{v_r}_{L^\infty((0,R)\times(0,T_0))}\leq CR<\infty,
\]
verifying the first part of the lemma.\\
If moreover $n\geq 2$, then the second equation in (\ref{DSR}) yields
\[
(r^{n-1}v_r)_r=r^{n-1}\mu+r^{n-1}u
\]
and thus upon integration for $r\in (0,R]$ and $\delta\in(0,r)$
\[
r^{n-1}v_r(r,t)-\delta^{n-1}v_r(\delta,t)=\f{\mu r^n}n-\f{\mu \delta^n}n-\int\limits_\delta^r \rho^{n-1}u(\rho,t)d\rho.
\]
Since $n-1>0$ and $v_r \in L^\infty((0,R)\times (0,T_0))$, taking $\delta \searrow 0$ results in
\[
r^{n-1}v_r(r,t)=\f{\mu r^n}n-\int\limits_0^r \rho^{n-1}u(\rho,t)d\rho,
\]
which upon dividing both sides by $r^{n-1}$ gives rise to (\ref{vr}). \qed
\vspace{20pt}

Under a weak assumption on $\beta$, classical solutions to $(\ref{DSR})$ conserve $\norm {u(\cdot,t)} _{L^1((0,R))}$ for at least as long as $u$ is bounded and (\ref{vr}) holds. These additional conditions are necessary since in contrast to usual settings our solution is not defined on a compact space.\\

\begin{lem}\label{lemmassc}
Suppose that $n\geq 1,~R>0$, $\beta>2-n$, and $u_0$ fulfills (\ref{init}). Let $(u,v)$ be a classical solution of (\ref{DSR}) in the sense of Definition \ref{classsol}, and assume that additionally 
\begin{equation}
0\leq u\in L^\infty((0,R)\times (0,T_0)) \qquad \text{for some }T_0\in(0,T] \label{ubound}
\end{equation}
as well as (\ref{vr}) holds. Then the mass conservation property
\be{massc}
\int_0^R \rho^{n-1} u(\rho,t) d\rho = \int_0^R \rho^{n-1} u_0(\rho) d\rho
\ee
is valid for all $t\in[0,T_0)$.
\end{lem}

\Proof
Let $(\zeta^{(\delta)})_{\delta\in(0,\f R2)}$ be a family of cutoff functions such that for all $\delta\in(0,\f R2)$ we have that $\zeta^{(\delta)} \in C^\infty([0,R])$ satisfies
\begin{equation}\label{zeta}
\begin{cases}
\zeta^{(\delta)}(r)=0, \qquad &r\in[0,\f\delta 2],\\
0\leq \zeta^{(\delta)}(r) \leq 1, \qquad &r\in(\f\delta 2,\delta),\\
\zeta^{(\delta)}(r)=1, \qquad &r\in[\delta,R],
\end{cases}
\end{equation}
as well as
\begin{equation}
0 \leq \zeta^{(\delta)}_r(r) \leq \f4\delta, \qquad r\in\bigg(\f\delta 2,\delta\bigg),
\end{equation}
and for some $C_1>0$ independent of $\delta$
\begin{equation}
|\zeta^{(\delta)}_{rr}(r)| \leq \f {C_1}{\delta^2}, \qquad r\in\bigg(\f\delta 2,\delta\bigg).
\end{equation}
Since (\ref{zeta}) guarantees that for all $\delta\in(0,\f R2)$ and $t\in(0,T_0)$ we have $\zeta^{(\delta)} u_t(\cdot,t) \in L^1((0,R))$, using (\ref{DSR}) we may compute
\begin{align}\label{massdiff}
\f{d}{dt}\int_0^R r^{n-1}\zeta^{(\delta)} u dr &= \int_0^R r^{n-1} \zeta^{(\delta)} u_t dr \notag\\
&=\int_0^R \zeta^{(\delta)} \cdot (r^{n-1+\beta}u_r -r^{n-1} uv_r)_r dr \notag\\
&=\int_0^R (\zeta^{(\delta)}_r \cdot r^{n-1+\beta})_r u dr + \int_0^R \zeta^{(\delta)}_r \cdot r^{n-1} uv_r dr
\end{align}
via partial integration. Herein, abbreviating $C_2:=\norm{u}_{L^\infty((0,R)\times(0,T_0))}$, we further estimate
\begin{align}\label{1stsumm}
\bigg|\int_0^R (\zeta^{(\delta)}_r \cdot r^{n-1+\beta})_r u dr\bigg|&\leq (n-1+\beta) \int_\f\delta 2^\delta |r^{n-2+\beta}\zeta^{(\delta)}_r u|+\int_\f\delta 2^\delta |r^{n-1+\beta}\zeta^{(\delta)}_{rr} u| \notag\\
&\leq (n-1+\beta) C_2 \int_\f\delta 2^\delta |r^{n-2+\beta}\zeta^{(\delta)}_r|+ C_2 \int_\f\delta 2^\delta |r^{n-1+\beta}\zeta^{(\delta)}_{rr}| \notag\\
&\leq 2(n-1+\beta) C_2 \delta^{n-2+\beta}+ C_2 C_1 \cdot \f12 \delta^{n-2+\beta},
\end{align}
for $t\in(0,T_0)$ by $n-2+\beta>0$, also guaranteeing the right hand side converges towards $0$ as we let $\delta \searrow 0$.\\
Since due to (\ref{vr}), Lemma \ref{vrlem} ensures that for some $C_3>0$
\[
|v_r(r,t)|\leq C_3r \qquad \text{for all }(r,t)\in(0,R]\times(0,T_0),
\]
we deduce that
\begin{align}\label{2ndsumm}
\bigg|\int_0^R \zeta^{(\delta)}_r \cdot r^{n-1} uv_r dr\bigg| &\leq C_2C_3\int_0^R \zeta^{(\delta)}_r \cdot r^n dr \notag\\
&\leq \f{4C_2C_3}n \cdot \delta^n \notag\\
&\longrightarrow 0,
\end{align}
for $\delta \searrow 0$, since $n\geq 1$. For $T\in(0,T_0)$, integrating (\ref{massdiff}) over $(0,T)$ now yields
\begin{align*}
\int_0^R r^{n-1}\zeta^{(\delta)} u(r,T) dr - \int_0^R r^{n-1}\zeta^{(\delta)} u_0(r) dr =& \int_0^T\int_0^R (\zeta^{(\delta)}_r \cdot r^{n-1+\beta})_r u drdt \\&+ \int_0^T\int_0^R \zeta^{(\delta)}_r \cdot r^{n-1} uv_r drdt.
\end{align*}
By (\ref{1stsumm}), (\ref{2ndsumm}) and monotone as well as dominated convergence, letting $\delta \searrow 0$ this results in
\[
\int_0^R r^{n-1} u(r,T) dr - \int_0^R r^{n-1}u_0(r) dr = 0,
\]
verifying (\ref{massc}). \qed
\abs

Of major importance to our further analysis is the transformation of (\ref{DSR}) to a Dirichlet problem.

\begin{lem}\label{dirtraf}
Suppose $n\geq 1,~R>0$, $\beta>2-n$, and $u_0$ fulfills (\ref{init}), and let $(u,v)$ be a classical solution to (\ref{DSR}) in $(0,R]\times[0,T)$ for which (\ref{vr}) and (\ref{ubound}) hold for $T_0=T$.\\
We introduce the mass accumulation function
\begin{equation}
w(s,t):=\int\limits_0^{s^\frac{1}{n}} \rho^{n-1}u(\rho,t)d\rho,~~s=r^n \in [0,R^n], t\in [0,T). \label{w}
\end{equation}
Then $w \in C^0([0,R^n]\times[0,T))\cap C^{2,1}((0,R^n]\times(0,T))$, and for all $s\in (0,R^n)$ and $t\in (0,T)$, its spatial derivatives are given by
\begin{equation}
w_s(s,t)=\frac{1}{n}\cdot u(s^\frac{1}{n},t)~~~\text{and}~~~w_{ss}(s,t)=\frac{1}{n^2}\cdot s^{\frac{1}{n}-1}u_r(s^{\frac{1}{n}},t). \label{ws}
\end{equation}
Furthermore, $w$ solves
the Dirichlet problem 
\be{0w}
    	\left\{ \begin{array}{ll}
	w_t = n^2 s^{2-\frac{2}{n}+\f{\beta}n} w_{ss} + n ww_s - \mu s w_s,
	\qquad & s\in (0,R^n), \ t\in (0,T), \\[1mm]
	w(0,t)=0, \quad w(R^n,t)=\frac{m}{\omega_n},
	\qquad & t\in (0,T), \\[1mm]
	w(s,0)=w_0(s),
	& s\in (0,R^n),
 	\end{array} \right.
\ee
with $m:=\io u_0$, $\mu=\frac{nm}{\omega_n R^n}$ and
\be{w0}
	w_0(s):=\int_0^{s^\frac{1}{n}} \rho^{n-1} u_0(\rho) d\rho,
	\qquad s\in [0,R^n].
\ee
\end{lem}

\Proof
Since $u$ and $v_r$ are bounded and mass conservation holds, this can be obtained via straightforward calculation as in \cite[p.264]{Win2011}.\qed

\newpage

\subsection{Constructing a solution to the scalar Dirichlet problem}\label{auxproblem}

By Lemma \ref{dirtraf}, under the mild assumption that $\beta>2-n$, every bounded nonnegative classical solution of (\ref{DSR}) satisfying (\ref{vr}) implies the existence of a classical solution to (\ref{0w}). Thus, in search of solutions to the former, it appears sensible to study the latter system.\\
Regarding similar systems emerging from the basic parabolic-elliptic Keller-Segel model in two dimensions, proofs for the existence of corresponding solutions have been sketched in \cite{Win2019} and \cite{biler08}.\\

As (\ref{0w}) still contains a diffusion degeneracy, questions regarding its solvability are not covered by standard theory. Therefore, we resort to non-degenerate \glqq approximating\grqq~problems. Their form is chosen so as to allow for concavity later.  For $\eps \in (0,R^n)$, consider
\be{weps}
    	\left\{ \begin{array}{ll}
	w_{\eps t} = n^2 s^{2-\frac{2}{n}+\f{\beta}n} w_{\eps ss} + n w_\eps w_{\eps s} - \mu (s-\eps) w_{\eps s},
	\qquad & s\in (\eps,R^n), \ t>0, \\[1mm]
	w_\eps(\eps,t)=0, \quad w_\eps(R^n,t)=\frac{m}{\omega_n},
	\qquad & t>0, \\[1mm]
	w_\eps(s,0)=w_{0\eps}(s),
	& s\in (\eps,R^n).
 	\end{array} \right.
\ee
Here and in the following, we always let $m>0$, which can be interpreted as the total mass in the original system as in Lemma \ref{dirtraf}.\\

We can obtain a solution to -- albeit yet an incomplete version of -- the degenerate problem formulated in (\ref{0w}):

\begin{lem}\label{limitfunction}
Let $n\geq 1,~R>0$, $\theta\in(0,1)$, $\mu=\frac{nm}{\omega_n R^n}$, $\beta>0$ and suppose $w_0 \in C^{1+\theta}([0,R^n])$ is such that 
\begin{equation}\label{w0def}
w_0(0)=0, \qquad w_0(R^n)=\f m{\omega_n} \quad \text{as well as} \quad w_{0s}\geq 0.
\end{equation}
Then there exists a function $w\in C^{1,\f12}((0,R^n]\times [0,\infty)) \cap C^{2,1}((0,R^n]\times(0,\infty))$ which solves the incomplete Dirichlet problem
\be{inc0w}
    	\left\{ \begin{array}{ll}
	w_t = n^2 s^{2-\frac{2}{n}+\f{\beta}n} w_{ss} + n ww_s - \mu s w_s,
	\qquad & s\in (0,R^n), \ t>0, \\[1mm]
	w(R^n,t)=\frac{m}{\omega_n},
	\qquad & t>0, \\[1mm]
	w(s,0)=w_0(s),
	& s\in (0,R^n).
 	\end{array} \right.
\ee
Moreover, $w$ is bounded via
\[
0 \leq w(s,t) \leq \f m{\omega_n} \qquad \text{for all}~(s,t) \in (0,R^n]\times[0,\infty),
\]
and $w(\cdot,t)$ is monotonically increasing in $(0,R^n]$ for all $t \in[0,\infty)$. Thus, for fixed $t\geq 0$ we may continuously extend $w(\cdot,t)$ to $s=0$ via
\be{wext}
w(0,t):=\lim_{s \searrow 0} w(s,t) \geq 0.
\ee
\end{lem}

\Proof
For $\eps\in(0,R^n)$, define $w_{0\eps} \in C^{1+\theta}([\eps,R^n])$ as
\begin{equation}\label{w0epsdef}
w_{0\eps}(s):=w_0\Big(\f{R^n(s-\eps)}{R^n-\eps}\Big), \qquad s\in[\eps,R^n].
\end{equation}
This warrants that $w_{0\eps}(\eps)=0,~w_{0\eps}(R^n)=\f m{\omega_n},~w_{0\eps s}\geq 0$ in $[\eps,R^n]$ as well as $w_{0\eps}\leq w_0$ in $[\eps,R^n]$, $w_{0\eps} \nearrow w_0$ in $(0,R^n]$ and $w_{0\eps} \to w_0$ in $C^1_{loc}((0,R^n])$ as $\eps \searrow 0$.
Then the system (\ref{weps}) possesses a unique global classical solution $w_\eps \in C^0([\eps,R^n]\times[0,\infty))\cap C^{2,1}([\eps,R^n]\times(0,\infty))$ (cf. \cite[V.6]{LSU}).\\
 Via two applications of (\ref{comparison}), $w_\eps\leq \f m{\omega_n}$ and $(w_\eps)_{\eps\in (0,R^n)}$ is pointwise monotically decreasing in $\eps$, so these solutions satisfy $w_\eps \nearrow w$
  in $(0,R^n] \times [0,\infty)$ as $\eps\searrow 0$ with some limit function $w$
  fulfilling $ 0 \le w \le \frac{m}{\omega_n}$, since $0 \le w_\eps \le \f m {\omega_n}$ 
  in $[\eps,R^n] \times [0,\infty)$ by another evident comparison argument.
  Because $w_{0\eps s} \ge 0$ and $w_{\eps s} \ge 0$ both on $\{\eps\} \times (0,\infty)$
  and on $\{R^n\} \times (0,\infty)$, we may furthermore conculde that $w_{\eps s} \ge 0$ in $(\eps,R^n) \times (0,\infty)$.
  From interior parabolic Schauder estimates (\cite{LSU}) and the Arzel\`a-Ascoli theorem, we furthermore obtain
  that in fact $w_\eps \to w$ in $C^{1,\f12}_{loc} ((0,R^n] \times [0,\infty))$ and in $C^{2,1}_{loc}((0,R^n] \times (0,\infty))$
  as $\eps\searrow 0$, and that it solves (\ref{inc0w}) classically. \qed
\abs

Note that $w$ is continuous at $(0,t)$ with respect to the spatial variable $s$ for fixed time $t$, yet not necessarily continuous in spacetime. As results in related systems, such as \cite[Theorem 3.1]{biler08} or \cite[Lemma 3.4]{primzaltz}, suggest, in this general framework this is probably not always the case globally in time. For small times however, we can establish

\begin{lem}\label{w0bound}
Let $n\geq 1,~R>0$, $\theta\in(0,1)$, $\mu=\frac{nm}{\omega_n R^n}$, $\beta>0$ and suppose $w_0 \in C^{1+\theta}([0,R^n])$ is as in (\ref{w0def}).\\
Then there is a $T^*>0$ and $y:[0,T^*) \to \R$ such that with $w$ as in Lemma \ref{limitfunction}
\be{compfunc}
w(s,t) \leq y(t) \cdot s \qquad \text{for all} \quad (s,t)\in (0,R^n] \times [0,T^*).
\ee
In consequence, we obtain that
\be{w00}
w\in C^0([0,R^n]\times[0,T^*))\quad \text{with} \quad w(0,t)=0 \qquad \text{for all} \quad t \in [0,T^*).
\ee
\end{lem}

\Proof
Once more, we resort to a comparison argument utilizing the approximating solutions $w_\eps$ which shall again be given by (\ref{w0epsdef}) and (\ref{weps}).\\
Define 
\be{initialODE}
y_0:=\|w_{0s}\|_{L^\infty((0,R^n))},
\ee
and let $y\in C^1([0,T^*))$ denote the solution to the ODE system
\be{odey}
    	\left\{ \begin{array}{ll}
	y'(t) = ny^2(t),
	\qquad & t>0, \\[1mm]
	y(0)=y_0,
 	\end{array} \right.
\ee
extended up to its maximal time of existence $T^*:=T^y_{max}>0$.\\
Based thereupon, for $\eps \in (0,R^n)$ and arbitrary $T\in (0,T^*)$, set 
\be{ow}
\ow(s,t):=y(t) \cdot s, \qquad (s,t)\in [\eps,R^n]\times[0,T].
\ee
Then we have
\begin{align} \label{diffit}
\ow_t - n^2 s^{2-\frac{2}{n}-\f{\beta}n} \ow_{ss} - n \ow \ow_s + \mu (s-\eps) \ow_s &= y_t\cdot s - n y^2 \cdot s + \mu y \cdot (s-\eps) \notag\\
&= ny^2 \cdot s - ny^2 \cdot s + \mu y \cdot (s-\eps) \notag \\
&= \mu y \cdot (s-\eps)\\
&\geq 0
\end{align}
in $(\eps,R^n)\times (0,T)$. Furthermore, by (\ref{w0def}), (\ref{w0epsdef}) and (\ref{initialODE}) we deduce that for $s \in (\eps,R^n)$
\begin{align}\label{w0epslower}
w_\eps(s,0)&=w_{0\eps}(s) \leq w_0(s) \notag\\
&= \underbrace{w_0(0)}_{=0} + \int_0^s w_{0s}(\rho) d\rho \notag\\
&\leq \|w_{0s}\|_{L^\infty((0,R^n))} \cdot s \notag\\
&\leq y_0 \cdot s \notag\\
&= \ow(s,0).
\end{align}
Since $y_0\geq \f{m}{\omega_n R^n}= \f{\mu}n$, by means of a simple ODE comparison argument we can see that $y(t)\geq \f{\mu}n \geq 0$ for all $t\geq 0$.
Therefore, the required inequalities at the lateral boundary are easily confirmed via
\be{epsbd}
w_\eps(0,t)=0 \leq y(t) \cdot \eps = \ow(\eps,t)
\ee
and
\be{Rnbd}
w_\eps(R^n,t)= \f m{\omega_n}= \f{\mu}n \cdot R^n \leq y(t) \cdot R^n = \ow(R^n,t)
\ee
for $t \in (0,T)$.\\
Thus combining (\ref{diffit}), (\ref{w0epslower}), (\ref{epsbd}) and (\ref{Rnbd}), an application of Lemma \ref{comparison} yields
\[
w_\eps (s,t) \leq \ow(s,t) = y(t) \cdot s \qquad \text{for all} \quad (s,t)\in [\eps,R^n]\times [0,T],
\]
and by taking $T \nearrow T^*$ this inequality holds in $[\eps,R^n]\times [0,T^*)$. Since the right hand side is independent of $\eps$ and $\eps \in (0,R^n)$ has been chosen arbitrarily, this results in (\ref{compfunc}).\\
This obviously entails (\ref{w00}) as well.\qed
\abs

Note that Lemma \ref{w0bound} implies that $w$ solves (\ref{0w}) in $[0,R^n]\times[0,T^*)$.

The information on the boundary $s=0$ on hand now allows for a statement on uniqueness of the function constructed in Lemma \ref{limitfunction} as a solution to (\ref{0w}), as long as its spatial derivative is bounded.

\begin{lem}\label{uniqueness}
Let $n\geq 1,~R>0$, $\theta\in(0,1)$, $\mu=\frac{nm}{\omega_n R^n}$, $\beta>0$ and suppose that $w_0 \in C^{1+\theta}([0,R^n])$ is as in (\ref{w0def}).\\
For the function $w$ defined in Lemma \ref{limitfunction}, let $T>0$ be such that $w \in C^0([0,R^n]\times[0,T))$ with 
\be{w000}
w(0,t)=0 \qquad \text{for all} \quad t\in(0,T),
\ee
and furthermore assume that $w_s$ is bounded in $(0,R^n)\times[0,T)$.\\
Then $w$ is the unique solution of (\ref{0w}) in $C^0([0,R^n]\times[0,T))\cap C^{2,1}((0,R^n]\times(0,T))$.
\end{lem}

\Proof
By (\ref{inc0w}) and (\ref{w000}), $w$ is a classical solution of (\ref{0w}) in $[0,R^n]\times[0,T)$.\\
Combined with the boundedness of $w_s$, via Lemma \ref{comparison} we may immediately infer uniqueness. \qed
\abs

\subsection{Local boundedness of $w_s$}\label{wssection}

The linear supersolution of $w$ established in Lemma \ref{w0bound} alone is not sufficient to deduce that the spatial derivative $w_s$ is bounded in $(0,R^n)\times(0,T^*)$. In order to be able to draw this conclusion, we impose stricter requirements on the initial data to ensure $w(\cdot,t)$ is concave.

\begin{lem}\label{concavitytw}
Let $n\geq 1,~R>0$, $\mu=\frac{nm}{\omega_n R^n}$, $\theta \in (0,1)$, $\beta>0$ and $w_0 \in C^{2+\theta}([0,R^n])$ be as in (\ref{w0def}) with
\be{w0ssneg}
w_{0ss}(s) \leq 0, \qquad s \in (0,R^n),
\ee
as well as
\be{compatibility}
w_{0ss}(0)=0, \quad w_{0s}(R^n)=0 \quad \text{and} \quad w_{0ss}(R^n)=0.
\ee
In that case, $w \in C^{2,1}((0,R^n]\times[0,\infty))$ as in Lemma \ref{limitfunction} satisfies
\be{twssle0}
w_{ss}(s,t) \leq 0 \qquad \text{for all} \quad (s,t) \in (0,R^n)\times [0,\infty).
\ee
\end{lem}

\Proof
For $\eps\in(0,R^n)$, define $w_{0\eps} \in C^{1+\theta}([\eps,R^n])$ via (\ref{w0epsdef}) and let $w_\eps \in C^0([\eps,R^n]\times[0,\infty))\cap C^{2,1}([\eps,R^n]\times(0,\infty))$.\\
The key to deducing concavity is ensuring that $w_\eps \in C^{2,1}([\eps,R^n]\times [0,\infty))$ to control the boundary near $t=0$.
To that end, by (\ref{compatibility}) we may infer that
\[
w_{0\eps ss}(\eps)=0, \quad w_{0\eps s}(R^n)=0 \quad \text{and} \quad w_{0\eps ss}(R^n)=0
\]
and thus
\[
n^2 \eps^{2-\frac{2}{n}+\f{\beta}n} \tilde w_{0\eps ss}(\eps) + n \tilde w_{0\eps}(\eps) \tilde w_{0\eps s}(\eps) - \mu (\eps-\eps) \tilde w_{0\eps s}(\eps)=0
\]
as well as
\[
n^2 R^{n(2-\frac{2}{n}+\f{\beta}n)} \tilde w_{0\eps ss}(R^n) + n \tilde w_{0\eps}(R^n) \tilde w_{0\eps s}(R^n) - \mu (R^n-\eps) \tilde w_{0\eps s}(R^n)=0.
\]
Therefore, the compatibility conditions are satified and thus the desired regularity can be derived via Schauder estimates.\\
Since (\ref{w0ssneg}) warrants that also $w_{0\eps ss}\leq 0$, a standard comparison argument yields
\be{wepsssle0}
w_{\eps ss}(s,t) \leq 0 \qquad \text{for all} \quad (s,t) \in (\eps,R^n)\times [0,\infty),
\ee
whence we may readily infer (\ref{twssle0}). \qed
\abs

\begin{lem}\label{concavityw}
Let $n\geq 1,~R>0$, $\mu=\frac{nm}{\omega_n R^n}$, $\theta \in (0,1)$, $\beta>0$ and $w_0 \in C^{2+\theta}([0,R^n])$ be as in (\ref{w0def}) with (\ref{w0ssneg}) as well as (\ref{compatibility}).\\
Let $w$ denote the global solution to (\ref{inc0w}) from Lemma \ref{limitfunction}. Then with $T^*>0$ as in Lemma \ref{w0bound} we have that for each $T\in(0,T^*)$ there exists $C=C(T)$ such that
\be{wslocbound}
w_s(s,t)\leq C \qquad \text{for all} \quad (s,t)\in (0,R^n]\times[0,T].
\ee
\end{lem}

\Proof
Since the conditions of Lemma \ref{concavitytw} are met, $w \in C^{2,1}((0,R^n]\times[0,\infty))$ fulfills
\be{wssle0}
w_{ss}(s,t) \leq 0 \qquad \text{for all} \quad (s,t) \in (0,R^n)\times [0,\infty).
\ee
Now let $0<T<T^*$. Lemma \ref{w0bound} ensures that
\be{upperline}
w(s,t)\leq C\cdot s \qquad \text{for all} \quad (s,t)\in(0,R^n]\times[0,T]
\ee
for some $C=C(T)>0$. Thus necessarily for $t\in(0,T)$
\[
\liminf_{s \searrow 0}w_s(s,t)\leq C,
\]
and seeing that (\ref{wssle0}) warrants $w_s(\cdot,t)$ to be monotonically decreasing for all $t\in[0,\infty)$, we may easily infer (\ref{wslocbound}). \qed

\subsection{Retransformation to Keller-Segel type system}\label{re}

With a local-in-time solution of (\ref{0w}) on hand, we may now obtain a solution to our original problem (\ref{DS}), or rather to its counterpart (\ref{DSR}) in radial coordinates.\\

\begin{lem}\label{retransform}
Let $n\geq 1,~R>0$, $\theta\in(0,1)$, $\mu=\frac{nm}{\omega_n R^n}$, $\beta>0$ and suppose that $w_0 \in C^{1+\theta}([0,R^n])$ is as in (\ref{w0def}).\\
Furthermore, we choose $T_0>0$ maximally such that the function 
\[w\in C^{1,\f12}((0,R^n]\times [0,\infty)) \cap C^{2,1}((0,R^n]\times(0,\infty))
\]
constructed in Lemma \ref{limitfunction} is in $C^0([0,R^n]\times[0,T_0))$ with
\be{leftiezero}
w(0,t)=0 \qquad \text{for all} \quad t\in(0,T_0).
\ee
Then for $u_0 \in C^\theta([0,R])$ defined via
\[
u_0(r)=n \cdot w_{0s}(r^n), \qquad r\in[0,R],
\]
the pair of functions $u \in C^0((0,R] \times [0,T_0))\cap C^{2,1}((0,R] \times (0,T_0))$ given by
\be{defu}
u(r,t)=n \cdot w_s(r^n,t),
\ee
and $v \in C^{2,0}((0,R] \times (0,T_0))$ fulfilling
\begin{equation}\label{vr1}
v_r(r,t)=\f1{r^{n-1}}\bigg(\f{\mu r^n}n-\int\limits_0^r \rho^{n-1}u(\rho,t)d\rho\bigg), \quad (r,t)\in (0,R] \times (0,T_0),
\end{equation}
for all $t\in(0,T_0)$ solves
\be{DSR1}
    	\left\{ \begin{array}{rcll}
	u_t &=& \f1{r^{n-1}}(r^{n-1+\beta}u_r)_r - \f1{r^{n-1}}(r^{n-1} uv_r)_r, 
	\qquad & r\in(0,R), \ t>0, \\[1mm]
	0 &=& \f1{r^{n-1}}(r^{n-1}v_r)_r - \mu + u,
	\qquad & r\in(0,R), \ t>0, \\[1mm]
	& & \hspace*{-10mm}
	u_r=v_r=0,
	\qquad & r=R, \ t>0, \\[1mm]
	& & \hspace*{-10mm}
	u(r,0)=u_0(r),
	\qquad & r\in(0,R),
 	\end{array} \right.
\ee
classically in $(0,R]\times[0,T_0)$.\\
Moreover, with $T^*>0$ as in Lemma \ref{w0bound},
\be{tmget*}
T_0 \geq T^*.
\ee
Furthermore, $u$ is nonnegative, and the total mass is conserved, that is
\be{massconserve}
\int_0^R \rho^{n-1} u(\rho,t) d\rho = \int_0^R \rho^{n-1} u_0(\rho) d\rho
\ee
for all $t\in(0,T_0)$.
\end{lem}

\Proof
Since $w \in C^{2,1}((0,R^n]\times(0,\infty))$ and
\[
w_t = n^2 s^{2-\f2n+\f\beta n} w_{ss} + n w w_s - \mu s w_s
\]
in $(0,R^n)\times (0,\infty)$, for any $t>0$ we necessarily have
\[
w_t(R^n,t)=n^2 R^{n(2-\f2n+\f\beta n)} w_{ss}(R^n,t) + n w(R^n,t) w_s(R^n,t) - \mu R^n w_s(R^n,t).
\]
Due to $w(R^n,t)=\f m{\omega_n}=\f{\mu R^n}n$ for all $t>0$ and thus also $w_t(R^n,t)=0$ in $(0,\infty)$, this yields
\[
0=n^2 R^{n(2-\f2n+\f\beta n)} w_{ss}(R^n,t) + w_s(R^n,t)\cdot \bigg(n \f{\mu R^n}n - \mu R^n\bigg),
\]
and therefore
\be{wss}
w_{ss}(R^n,t)=0 \qquad \forall t>0.
\ee
This yields the boundary condition
\be{ur}
u_r(R,t)=n^2 w_{ss}(R^n,t)\cdot R^{n-1}=0 \qquad \text{for} \quad t\in(0,T_0).
\ee
The remaining part can be verified by the asserted regularity and straightforward calculation. \qed
\abs

Combining multiple results of this section, we are moreover able to formulate a proposition regarding uniqueness of solutions to (\ref{DSR1}).

\begin{lem}\label{uni}
Let $n\geq 2$, $R>0$, $\theta\in(0,1)$, $\mu=\f{nm}{\omega_n R^n}$, $\beta>0$ and $u_0 \in C^\theta([0,R])$.\\
Then for $T>0$ there is at most one solution $(u,v)$ of (\ref{DSR1}) in $(0,R]\times[0,T)$ with
\[
\begin{cases}
u \in C^0((0,R] \times [0,T))\cap C^{2,1}((0,R] \times (0,T)),\\
v \in C^{2,0}((0,R] \times (0,T)),
\end{cases}
\]
which has the properties that $\int_0^R v(r,t) dr = 0$ for all $t\in(0,T)$ and
\be{unicond}
0\leq u \in L^\infty((0,R)\times (0,T)) \quad \text{and} \quad v_r\in L^\infty((0,R)\times(0,T)).
\ee
\end{lem}

\Proof
Note that $n\geq 2$ implies that $\beta>0\geq 2-n$. Therefore, if also (\ref{unicond}) holds, the requirements of Lemma \ref{vrlem} and Lemma \ref{dirtraf} are met and thus the existence of a solution of (\ref{DSR1}) implies the existence of a solution $w \in C^0([0,R^n]\times[0,T))\cap C^{2,1}((0,R^n]\times(0,T))$ of (\ref{0w}) with $w_0\in C^{1+\theta}([0,R^n])$, in which $w_s\in C^0([0,R^n]\times[0,T))$ is nonnegative and bounded.\\
Then Lemma \ref{uniqueness} however warrants that $w$ is the unique classical solution of (\ref{0w}). Via the fundamental theorem of calculus, we can also easily infer that then
\[
w(s,t)=w(0,t)+\int_0^s w_s(\rho,t) d\rho \leq \|w_s\|_{L^\infty((0,R^n)\times(0,T)}\cdot s,
\]
providing a pendant to (\ref{compfunc}), although we technically do not need it since we do not demand $v_r$ to be extendable to $r=0$.\\
Lemma \ref{retransform} then guarantees the existence of a solution $(u,v)$ of (\ref{DSR1}) in the sense specified above.\\
In conclusion, solutions of (\ref{DSR1}) satisfying the conditions of this lemma correspond with solutions of (\ref{0w}), thus transfering their uniqueness. \qed
\abs

We remark that actually $u_0 \in C^0([0,R])$ is sufficient for Lemma \ref{uni} to hold, since bounded $w_s \in C^0((0,R^n]\times[0,T))$ is sufficient for all relevant arguments.\\

Sharpening the conditions in accordance with subsection \ref{wssection} enables us to acquire additional properties, most prominently local-in-time boundedness of $u$ in $(0,R]\times[0,T^*)$.

\begin{lem}\label{uhu}
Suppose that the conditions of Lemma \ref{concavityw} hold, and let $(u,v)$ denote the solution of (\ref{DSR1}) constructed in Lemma \ref{retransform}.\\
Then this solution has the additional properties that $u \in C^{1,0}((0,R] \times [0,T_0))$ and
\be{urle0}
u_r(r,t)\leq 0 \qquad \text{for all} \quad (r,t)\in(0,R]\times[0,T_0).
\ee
Moreover, for $T^*>0$ as in Lemma \ref{w0bound} and each $T\in(0,T^*)$, there exists $C=C(T)>0$ such that
\be{realubound}
u(r,t)\leq C \qquad \text{for all} \quad (r,t)\in(0,R]\times[0,T].
\ee
\end{lem}

\Proof
One only needs to adapt the results of Lemma \ref{concavityw}.\\
The first spatial derivative of $u$ is given by
\[
u_r(r,t)=n^2 w_{ss}(r^n,t)\cdot r^{n-1} \qquad \text{for} \quad (r,t)\in(0,R]\times(0,T_0).
\]
Therefore the claimed regularity follows from $w\in C^{2,1}((0,R^n]\times[0,T_0)$, whereas (\ref{urle0}) is a consequence of (\ref{wssle0}).\\
Lastly, as $u(r,t)=n \cdot w_s(r^n,t)$ for $(r,t)\in(0,R]\times(0,T_0)$, (\ref{wslocbound}) translates to (\ref{realubound}). \qed

\newpage

\section{Ruling out global boundedness in (\ref{DS}) for sufficiently concentrated initial data}\label{bu}

In usual settings of the Keller-Segel system and its variants, the occurence of blow-up at a finite time $T<\infty$ corresponds with the maximal time of existence $T_{max}$ equaling $T$. Since our classical solution concept however does not require $u$ to be defined continuously on a compact space, it is well possible that for a domain $\Omega  \subset \R^n$ 
\[
\limsup_{t \nearrow T}\|u(\cdot,t)\|_{L^\infty(\Om)}=\infty
\]
but $u \in C^0(\Om\times[0,T_0)) \cap C^{2,1}(\Om\times(0,T_0))$ for some $T_0>T$. Moreover considering we have no extensibility criterion at hand, we restrict ourselves to ruling out the existence of global bounded solutions under certain circumstances.

As in \cite[Lemma 3.3]{Win2019}, we may establish that given a condition corresponding to the initial mass in the corresponding original system (\ref{DS}) being sufficiently concentrated, there is no global solution of (\ref{0w}) for which $w_s$ is bounded locally in time in $(0,R^n)\times(0,\infty)$.

\begin{lem}\label{lem87}
  Let $n\geq 1$, $\beta>\max\{0,2-n\}$, $R>0$, $m_0>0$ and $m\ge m_0$.
  There exists $s_0=s_0(m_0,m,R,\beta) \in (0,R^n)$ such that if $w_0 \in C^{1}([0,R^n])$,
  $\mu=\f{nm}{\omega_n R^n}$, and furthermore
  \be{87.1}
	w_0(s_0) \geq \frac{m_0}{\omega_n},
  \ee
  then there is no global classical solution 
\[
w\in C^0([0,R^n]\times[0,\infty))\cap C^{2,1}((0,R^n]\times(0,\infty))
\]
of (\ref{0w}) with the property that for each $T>0$
\be{wsglobbound}
w_s \in L^\infty((0,R^n)\times(0,T)).
\ee
\end{lem}
\Proof
Due to $\beta>2-n$, it is possible to fix $\gamma \in (0,1)$ with the property that
\be{gammarestr}
\gamma \leq 1-\f2n+\f\beta n.
\ee
  We abbreviate
  \bea{87.2}
	& &c_1:=\frac{8\bigg(2-\f2n+\f\beta n-\gamma\bigg)^2 n^3}{3-\frac{4}{n}+\f{2\beta}n-\gamma},
	\qquad 
	c_2:=\frac{2n}{(3-\gamma)\omega_n^2}
	\qquad \mbox{and} \nn\\
	& &c_3:=\frac{3}{4\omega_n} \cdot \Big(\f1{1-\gamma}-\f1{2-\gamma}\Big).
  \eea

Observe that $\beta>2-n$ implies that $2-\frac{4}{n}+\f{2\beta}n > 0$.\\
  Therefore, given $R>0$, $m_0>0$ and $m>0$ we can fix $s_0=s_0(m_0,m,R)\in (0,\frac{R^n}{2})$ such that $s_1:=2s_0$ satisfies
  \be{87.3}
	s_1^{2-\frac{4}{n}+\f{2\beta}n} \le \frac{(1-\gamma)c_3^2}{n c_1} m_0^2
  \ee
  and 
  \be{87.4}
	s_1^2
	\le \frac{(1-\gamma)c_3^2}{n c_2} \cdot \frac{m_0^2 R^{2n}}{m^2},
  \ee
  and henceforth we assume that $w\in C^0([0,R^n]\times[0,\infty))\cap C^{2,1}((0,R^n]\times(0,\infty))$ is a global classical solution of (\ref{0w}) satisfying (\ref{wsglobbound}).
  For $\delta\in (0,\frac{s_1}{2})$, we then use (\ref{0w}) to compute
  \bas
	& & \hspace*{-10mm}
	\frac{d}{dt} \int\limits_\delta^{s_1}
	s^{-\gamma} (s_1-s) w(s,t) ds \\
	&=& n^2 \int\limits_\delta^{s_1} s^{2-\f2n+\f\beta n - \gamma} (s_1-s) w_{ss}(s,t) ds
	+ \frac{n}{2} \int\limits_\delta^{s_1} s^{-\gamma}(s_1-s) (w^2)_s (s,t) ds \\
	& &- \mu \int\limits_\delta^{s_1} s^{1-\gamma} (s_1-s) w_s(s,t) ds \\
	&=& n^2\bigg(2-\f2n+\f\beta n -\gamma\bigg)\bigg(1-\f2n+\f\beta n-\gamma\bigg) \int\limits_\delta^{s_1} s^{-\f2n+\f\beta n-\gamma} w(s,t) ds\\
	& & - 2 n^2\bigg(2-\f2n+\f\beta n-\gamma\bigg) \int\limits_\delta^{s_1} s^{1-\f2n+\f\beta n-\gamma} w(s,t) ds
	+\f n 2 \cdot \gamma \int_\delta^{s_1} s^{-\gamma-1}(s_1-s)w^2(s,t) ds\\
	& &+ \f n 2 \int\limits_\delta^{s_1} s^{-\gamma}w^2(s,t) ds + \mu (1-\gamma)\int\limits_\delta^{s_1} s^{-\gamma} (s_1-s)w(s,t) ds - \mu \int\limits_\delta^{s_1}s^{1-\gamma} w(s,t) ds\\
	& &+ n^2\bigg(2-\f2n+\f\beta n-\gamma\bigg) \delta^{1-\f2n+\f\beta n-\gamma}(s_1-\delta)w(\delta,t)- n^2 \delta^{2-\f2n+\f\beta n-\gamma}(s_1-\delta)w_s(\delta,t) \\
	& &+n^2s_1^{2-\f2n+\f\beta n -\gamma}w(s_1,t)-n^2\delta^{2-\f2n+\f\beta n -\gamma}w(\delta,t)\\
	& & - \f n 2 \delta^{-\gamma}(s_1-\delta)w^2(\delta,t)  + \mu \delta^{1-\gamma} (s_1-\delta)w(\delta,t)
	\qquad \mbox{for all } t>0
  \eas
  with $\mu=\frac{nm}{\omega_n R^n}$. Note that due to (\ref{gammarestr}), the first summand in the last equality is nonnegative. Neglecting some other nonnegative summands as well and integrating in time leads to
  \bea{87.5}
	\int_\delta^{s_1} s^{-\gamma} (s_1-s) w(s,t) ds
	&\ge& \int_\delta^{s_1} s^{-\gamma} (s_1-s) w_0(s) ds \nn\\
	& & - 2n^2\bigg(2-\f2n+\f\beta n-\gamma\bigg) \int_{0}^t \int_\delta^{s_1} s^{1-\f2n+\f\beta n-\gamma}w(s,\tau) dsd\tau \nn \\
	& &+ \frac{n}{2} \int_{0}^t \int_\delta^{s_1} s^{-\gamma} w^2(s,\tau) dsd\tau \nn
	 - \mu \int_{0}^t \int_\delta^{s_1} s^{1-\gamma} w(s,\tau) dsd\tau \nn\\
	& & - n^2\delta^{2-\f2n+\f\beta n-\gamma} (s_1-\delta) \int_{0}^t w_s(\delta,\tau) d\tau \nn \\
	& & - n^2 \delta^{2-\f2n+\f\beta n-\gamma} \int_{0}^t w(\delta,\tau) d\tau \nn\\
	& & - \frac{n}{2}\delta^{-\gamma} (s_1-\delta) \int_{0}^t w^2(\delta,\tau) d\tau
	\qquad \mbox{for all } t>0.
  \eea
Here since we assume $w(0,t)=0$ for all $t\geq 0$ and boundedness of $w_s$ in $(0,R^n) \times (0,t)$, we may infer that
  $\sup_{(s,\tau)\in (0,R^n)\times (0,t)} \frac{w(s,\tau)}{s}$ is finite, and therefore
  \bas
	& &n^2\delta^{2-\f2n+\f\beta n-\gamma} (s_1-\delta) \int_{0}^t w_s(\delta,\tau) d\tau
	+ n^2 \delta^{2-\f2n+\f\beta n-\gamma} \int_{0}^t w(\delta,\tau) d\tau \nn \\
	& &+ \frac{n}{2}\delta^{-\gamma} (s_1-\delta) \int_{0}^t w^2(\delta,\tau) d\tau
	\to 0
	\quad \mbox{as } \delta\searrow 0,
  \eas
  whence on several applications of the monotone convergence theorem we infer from (\ref{87.5}) that
  $y(t):=\int_0^{s_1} s^{-\gamma} (s_1-s) w(s,t) ds$, $t\ge 0$, satisfies
  \bas
	y(t)
	&\ge& y(0) 
	- 2n^2\bigg(2-\f2n+\f\beta n-\gamma\bigg) \int_{0}^t \int_0^{s_1} s^{1-\f2n+\f\beta n-\gamma}w(s,\tau) dsd\tau \\
	& &+ \frac{n}{2} \int_{0}^t \int_0^{s_1} s^{-\gamma}w^2(s,\tau) dsd\tau 
	- \mu \int_{0}^t \int_0^{s_1} s^{1-\gamma} w(s,\tau) dsd\tau
	\qquad \mbox{for all } t>0.
  \eas
  By Young's inequality,
  \bas
	& &2n^2\bigg(2-\f2n+\f\beta n-\gamma\bigg) \int_0^{s_1} s^{1-\f2n+\f\beta n-\gamma}w(s,\tau) ds\\
	&\le& \frac{n}{8} \int_0^{s_1}s^{-\gamma} w^2(s,\tau) ds
	+ 8\bigg(2-\f2n+\f\beta n-\gamma\bigg)^2 n^3 \int_0^{s_1} s^{2-\frac{4}{n}+\f{2\beta}n-\gamma} ds \\
	&=& \frac{n}{8} \int_0^{s_1} s^{-\gamma}w^2(s,\tau) ds
	+ c_1 s_1^{3-\frac{4}{n}+\f{2\beta}n-\gamma}
	\qquad \mbox{for all } \tau>0
  \eas
  and
  \bas
	\mu \int_0^{s_1} s^{1-\gamma} w(s,\tau) ds
	&\le& \frac{n}{8} \int_0^{s_1} s^{-\gamma} w^2(s,\tau) ds
	+ \frac{2\mu^2}{n} \int_0^{s_1} s^{2-\gamma} ds \\
	&=& \frac{n}{8} \int_0^{s_1} s^{-\gamma} w^2(s,\tau) ds
	+ c_2 \frac{m^2}{R^{2n}} s_1^{3-\gamma}
	\qquad \mbox{for all } \tau>0
  \eas
  as well as
  \bas
	y(\tau)
	&\le& \bigg\{ \int_0^{s_1} s^{-\gamma} w^2(s,\tau) ds \bigg\}^\frac{1}{2} \cdot
	\bigg\{ \int_0^{s_1} s^{-\gamma} (s_1-s)^2 ds \bigg\}^\frac{1}{2} \\
	&\le& \bigg\{ \int_0^{s_1} s^{-\gamma} w^2(s,\tau) ds \bigg\}^\frac{1}{2} \cdot
	\bigg\{ s_1^2 \int_0^{s_1} s^{-\gamma} ds \bigg\}^\frac{1}{2} \\
	&=& \bigg\{ \int_0^{s_1} s^{-\gamma}w^2(s,\tau) ds \bigg\}^\frac{1}{2} \cdot
	\Big\{ \frac{1}{1-\gamma} s_1^{3-\gamma} \Big\}^\frac{1}{2}
	\qquad \mbox{for all } \tau>0
  \eas
  by the Cauchy-Schwarz inequality. This entails that
  \be{87.6}
	y(t) \ge y(0)
	+ \frac{4(1-\gamma)}{n} s_1^{\gamma-3} \int_{0}^t y^2(\tau) d\tau
	- \Big\{ c_1 s_1^{3-\frac{4}{n}+\f{2\beta}n-\gamma} + c_2 \frac{m^2}{R^{2n}} s_1^{3-\gamma} \Big\} \cdot t
	\qquad \mbox{for all } t>0.
  \ee
  Now since (\ref{87.1}) along with our selections of $s_0$ and $c_3$ guarantees that
  \begin{align*}
	y(0)
	&\ge \frac{m_0}{\omega_n} \cdot \int_\frac{s_1}{2}^{s_1} s^{-\gamma} (s_1-s) ds \\
	&=\frac{m_0}{\omega_n} \cdot \bigg( \f1{1-\gamma}s_1\bigg(s_1^{1-\gamma}-\bigg(\f{s_1}2\bigg)^{1-\gamma}\bigg)-\f1{2-\gamma}			\bigg(s_1^{2-\gamma}-\bigg(\f{s_1}2\bigg)^{2-\gamma}\bigg)\bigg)\\
	&=\frac{m_0}{\omega_n} \cdot \bigg( \f3{4(1-\gamma)}s_1^{2-\gamma}-\f3{4(2-\gamma)}s_1^{2-\gamma}\bigg)\\
	&= c_3 m_0 s_1^{2-\gamma}
  \end{align*}
  and that hence, by (\ref{87.3}) and (\ref{87.4}),
  \bas
	\frac{c_1 s_1^{3-\frac{4}{n}+\f{2\beta}n-\gamma} + c_2 \frac{m^2}{R^{2n}} s_1^{3-\gamma}}{\f{2(1-\gamma)}{n} s_1^{\gamma-3} y^2(0)}
	\le \frac{n c_1}{2(1-\gamma)c_3^2 m_0^2} s_1^{2-\frac{4}{n}+\f{2\beta}n}
	+ \frac{n c_2 m^2}{2(1-\gamma)c_3^2 m_0^2 R^{2n}} s_1^2
	\le \frac{1}{2} + \frac{1}{2} =1,
  \eas
  it follows that there exists $T>0$ such that the problem
  \bas
	\left\{ \begin{array}{l}
	\uy'(t) = \frac{4(1-\gamma)}{n} s_1^{\gamma-3} \uy^2(t)
	- \Big\{ c_1 s_1^{3-\frac{4}{n}+\f{2\beta}n-\gamma} + c_2 \frac{m^2}{R^{2n}} s_1^{3-\gamma} \Big\}, 
	\qquad t\in (0,T), \\[1mm]
	\uy(0)=y(0),
	\end{array} \right.
  \eas
  admits a solution $\uy\in C^1([0,T))$ fulfilling $\uy(t)\nearrow +\infty$ as $t\nearrow T$.
  But an ODE comparison argument based on (\ref{87.6}) ensures that $y(t)\ge \uy(t)$ for all $t\in (0,T)$, which
  is incompatible with our hypothesis that $w$ is a global classical solution of (\ref{0w}) for which the first spatial derivative $w_s$ is bounded locally in 	time.
\qed

\newpage
\section{Proof of main results}

We shall now give proof to theorems \ref{maintheo1} -- \ref{maintheo3}.
\abs
The first two main theorems are obtained by utilizing the results from subsection \ref{re} and transfering them to (\ref{DS}).\\

{\sc Proof} of Theorem \ref{maintheo1}. \quad
Let $m:=\io u_0$. Consider $\tilde u_0 \in C^\theta([0,R])$ as the initial data in radial coordinates.
Then $w_0:[0,R^n]\to \R$ defined by 
\be{tilw}
w_0(s)=\int_0^{s^\f1n}\rho^{n-1}\tilde u_0 d\rho \qquad \text{for}\quad s\in[0,R^n]
\ee
is in $C^{1+\f \theta n}([0,R^n])$ and satisfies
\be{nw0s}
n \cdot w_{0s}(r^n)=\tilde u_0(r)
\ee
for $r\in[0,R]$ as well as 
\[
w_0(0)=0, \qquad w_0(R^n)=\f m{\omega_n} \quad \text{and} \quad w_{0s}(s)\geq 0, \quad s\in[0,R^n],
\]
and thus with $\mu:=\f{nm}{\omega_n R^n}$, Lemma \ref{retransform} asserts Theorem \ref{maintheo1}. \qed
\abs

For the second theorem, the key lemmata are \ref{uhu} and \ref{uni}.\\

{\sc Proof} of Theorem \ref{maintheo2}. \quad
Again, let $\tilde u_0 \in C^{1+\theta}((0,R])$ be the initial data in radial coordinates. Then $w_0$ as in (\ref{tilw}) is in $C^{2+\f \theta n}((0,R^n])$ and satisfies (\ref{nw0s}) as well as
\[
\tilde u_{0r}(r)=n^2 w_{0ss}(r^n)\cdot r^{n-1} \qquad \text{for} \quad r\in(0,R].
\]
Since $u_0$ is radially decreasing and therefore $\tilde u_{0r} \leq 0$ in $(0,R]$, this implies
\be{w0ssraddec}
w_{0ss}(s)\leq 0 \qquad \text{for all}\quad s\in(0,R^n].
\ee
Moreover, (\ref{Omboundary}) entails
\be{w0okay}
w_{0s}(R^n)=0 \qquad \text{and} \qquad w_{0ss}(R^n)=0.
\ee
Furthermore, due to $|\nabla u_0(x)| =|\tilde u_{0r}(r)|$ with $r=|x|$ for $x\in\bom$, (\ref{xlimit}) results in
\begin{align*}
|w_{0ss}(r^n)|&=\bigg|\f1{n^2}r^{1-n}\tilde u_{0r}(r)\bigg|\\
&\leq \f{C_0}{n^2}r^\theta \longrightarrow 0
\end{align*}
for $r\to 0$, warranting that $w_0\in C^{2+\f\theta n}([0,R^n])$ with
\be{w0annoying}
w_{0ss}(0)=0.
\ee
Now (\ref{w0ssraddec}), (\ref{w0okay}) and (\ref{w0annoying}) certify that the conditions of Lemma \ref{concavityw} are fulfilled, and thus Lemma \ref{uhu} warrants that $\tilde u \in C^{1,0}((0,R] \times [0,T_0))$ with
\be{turle0}
\tilde u_r(r,t)\leq 0 \qquad \text{for all} \quad (r,t)\in(0,R]\times[0,T_0),
\ee
and that moreover there exists $T^*\in(0,T_0]$ such that for $T\in(0,T^*)$,
\be{realtubound}
\tilde u(r,t)\leq C \qquad \text{for all} \quad (r,t)\in(0,R]\times[0,T]
\ee
with some $C=C(T)>0$. This however implies that for the solution $(u,v)$ of (\ref{DS}) from Theorem \ref{maintheo1}, (\ref{radiallydec}) and (\ref{realuboundtheo}) are valid.\\
By the definition of $\tilde v_r$ in (\ref{vr1}) and (\ref{realtubound}), we may infer that furthermore
\begin{align*}
|\tilde v_r(r,t)|&=\bigg|\f1{r^{n-1}}\bigg(\f{\mu r^n}n-\int\limits_0^r \rho^{n-1}\tilde u(\rho,t)d\rho\bigg)\bigg|\\
&\leq \f1{r^{n-1}}\cdot \f{\mu r^n}n + C(T)\f1{r^{n-1}}\int\limits_0^r \rho^{n-1} d\rho\\
&=\f{\mu r}n+C(T) \f r n
\end{align*}
for $T\in(0,T^*)$ and $(r,t)\in(0,R]\times(0,T)$, warranting
\be{vrboundproof}
\tilde v_r \in L^\infty((0,R)\times(0,T))
\ee
for all $T\in(0,T^*)$. This also entails that $\int_0^R \tilde v(r,t)dr$ is well-defined in $(0,T^*)$, and thus allows us to uniquely determine $\tilde v\in C^{2,0}((0,R] \times (0,T^*))$ by demanding
\[
\int_0^R \tilde v(r,t)dr=0.
\]
If now additionally $n\geq 2$, then by (\ref{realtubound}), (\ref{vrboundproof}) and nonnegativity of $\tilde u$, Lemma \ref{uni} ensures that $(\tilde u, \tilde v)$ is the unique solution of (\ref{DSR1}) in $(0,R]\times[0,T^*)$ satisfying these properties.\\
Regarding the original variables and considering
\[
|\nabla v(x,t)| =|\tilde v_{r}(|x|,t)| \quad \text{for} \quad (x,t)\in\bom\times(0,T^*)
\]
though, this means that the corresponding pair of functions $(u,v)$ is indeed the unique solution of (\ref{DS}) in $\Om_0\times[0,T^*)$ fulfilling (\ref{regularityclasses}) and (\ref{unicondtheo}).
\qed
\abs

The result on ruling out global boundedness is established by means of Subsection \ref{traforward} and Lemma \ref{lem87}.\\

{\sc Proof} of Theorem \ref{maintheo3}. \quad
Let $n\geq 2$, $R>0$, $\Om=B_R(0)\subset \R^n$, $\beta>0$ and $u_0\in C^0(\bom)$ complying with (\ref{init}). Then, with designations as above, $\tilde u_0 \in C^0([0,R])$.\\
Furthermore, suppose $u_0$ is such that with $m=\io u_0$, $m_0\in(0,m]$ and $r_0=s_0^\f1n$ for $s_0=s_0(m_0,m,R,\beta) \in (0,R^n)$ as in Lemma \ref{lem87}, we have that (\ref{initmasscon}) holds.\\
Assume there was a global classical solution to (\ref{DS}) satisfying (\ref{allisbound}) for all $T>0$ and (\ref{uvreginfty}).\\
Consider $(\tilde u, \tilde v)$ as $(u,v)$ in radial coordinates which has the property that
\[
\left\{ \begin{array}{l}
	\tilde u \in C^0((0,R] \times [0,\infty))\cap C^{2,1}((0,R] \times (0,\infty)),\\[1mm]
	\tilde v \in C^{2,0}((0,R] \times (0,\infty)),
	\end{array} \right.
\]
and solves (\ref{DSR}) for the initial condition $\tilde u(\cdot,0)=\tilde u_0$. Moreover, since $|\nabla v(x,t)|=|\tilde v_r(|x|,t)|$ for all $(x,t)\in\Om_0\times(0,\infty)$, (\ref{allisbound}) entails that for each $T>0$
\be{again}
0\leq \tilde u \in L^\infty((0,R)\times (0,T)) \quad \text{and} \quad \tilde v_r\in L^\infty((0,R)\times(0,T)).
\ee
As $n\geq 2$, Lemma \ref{vrlem} now ensures that $\tilde v_r$ is as in (\ref{vr}), whereby in turn we may infer that the requirements of Lemma \ref{dirtraf} are met, since $n\geq 2$ also implies that $\beta > 0 \geq 2-n$.\\
Therefore, $w:[0,R^n]\times[0,\infty)\to \R$ defined via
\[
w(s,t):=\int\limits_0^{s^\frac{1}{n}} \rho^{n-1}\tilde u(\rho,t)d\rho,~~s=r^n \in [0,R^n], \quad t\in [0,\infty),
\]
is in $C^0([0,R^n]\times[0,T))\cap C^{2,1}((0,R^n]\times(0,T))$ and solves (\ref{0w}) in $[0,R^n]\times[0,\infty)$ for $w_0\in C^1([0,R^n])$ given by
\[
	w_0(s)=\int_0^{s^\frac{1}{n}} \rho^{n-1} \tilde u_0(\rho) d\rho,
	\qquad s\in [0,R^n].
\]
Furthermore,
\[
w_s(s,t)=\frac{1}{n}\cdot \tilde u(s^\frac{1}{n},t) \qquad \text{for all}\quad (s,t)\in(0,R^n)\times(0,\infty)
\]
combined with (\ref{again}) entails that for each $T>0$
\[
w_s \in L^\infty((0,R^n)\times(0,T)).
\]
Thus a global classical solution of (\ref{0w}) exists, and its first spatial derivative is bounded locally in time. Since however for $s_0\in(0,R^n)$ as above by (\ref{initmasscon}),
\begin{align*}
w_0(s_0)&=\int_0^{s_0^\frac{1}{n}} \rho^{n-1} \tilde u_0(\rho) d\rho \\
&=\f1{\omega_n} \int_{B_{s_0^\f1n}(0)} u_0 \\
&=\f1{\omega_n} \int_{B_{r_0}(0)} u_0 \\
&\geq \f{m_0}{\omega_n},
\end{align*}
this is inconsistent with Lemma \ref{lem87}, thus ruling out the existence of a global classical solution to (\ref{DS}) satisfying (\ref{allisbound}) and (\ref{uvreginfty}) under the given premises. \qed

\newpage

\section{Appendix: A comparison principle for (\ref{0w})}
An important tool in the analysis of the transformed systems is a comparison principle. Since we also need to be able to deal with diffusion degeneracy, the standard comparison principles do not seem to apply. Also, for arguments as in \cite{Win2019}, some information drawn from standard results regarding the original Keller-Segel system is missing. Therefore, we prove a comparison principle tailored to our specific cases:
\abs

\begin{lem}\label{comparison}
Let $T>0$ and $l,L\geq 0$ with $l<L$. Suppose that $\underline w$ and $\overline w$ belong to $C^0([l,L]\times[0,T])\cap C^{2,1}((l,L)\times(0,T])$, and that additionally either 
\[
\uw_s \in L^\infty((l,L)\times(0,T)) \quad \text{or} \quad \ow_s \in L^\infty((l,L)\times(0,T)).
\]
Moreover, for $a,b,\gamma\geq 0$ and $\alpha,\delta,c,d\in\R$
\begin{equation}\label{dineq}
\uw_t\leq as^\alpha \uw_{ss}+bs^\gamma\uw\uw_s+cs^\delta \uw_s+d\uw_s \quad \text{and} \quad \ow_t\geq as^\alpha \ow_{ss}+bs^\gamma\ow\ow_s+cs^\delta \ow_s+d\ow_s
\end{equation}
shall hold for all $(s,t)\in (l,L)\times (0,T)$, as well as
\be{initineq}
\uw(s,0) \leq \ow(s,0) \qquad \text{for all }s\in(l,L)
\ee
and
\be{boundineq}
\uw(l,t) \leq \ow(l,t) \quad \text{and} \quad \uw(L,t) \leq \ow(L,t) \qquad \text{for all }t\in(0,T).
\ee
Then
\be{wineq}
\uw(s,t)\leq \ow(s,t) \qquad \text{for all } s\in[l,L] \quad \text{and} \quad t\in[0,T).
\ee
\end{lem}

\Proof
This can be derived via a standard argumentation to prove comparison principles for Dirichlet boundary problems. \qed

\vspace*{5mm}
{\bf Acknowledgement.} \quad
The author acknowledges support of the {\em Deutsche Forschungsgemeinschaft} in the context of the project
  {\em Fine structures in interpolation inequalities and application to parabolic problems}, project number 462888149.

\end{document}